\documentclass[a4,11pt]{article}

\usepackage{amsfonts}
\usepackage{amsmath}
\usepackage{amssymb}
\usepackage[mathscr]{eucal}
\usepackage{color}
\usepackage{hyperref}
\usepackage{amsthm,comment}
\theoremstyle{plain}
 \newtheorem{theorem}{Theorem}[section]
 
 \newtheorem*{theorem*}{Theorem}
 
 \newtheorem{lemma}[theorem]{Lemma}

\theoremstyle{definition}
 
\theoremstyle{remark}

\numberwithin{equation}{section}

\usepackage{mathtools}
\title{Uniqueness of an association scheme related to the Witt design on 11 points
}

\author{
Alexander L. Gavrilyuk
\and  
 Sho Suda
}

\date{\today}

\begin{document}

\maketitle

\abstract{
It follows from Delsarte theory that 
the Witt $4$-$(11,5,1)$ design gives rise to a $Q$-polynomial association 
scheme $\mathcal{W}$ defined on the set of its blocks. 
In this note we show that $\mathcal{W}$ is unique, i.e., 
defined up to isomorphism by its 
parameters.}

\section{Introduction}

A $t$-$(v,k,\lambda)$ design is a set of $k$-subsets (called {\bf blocks}) of $v$ points 
such that every $t$-subset is contained in exactly $\lambda$ blocks. 
A $t$-design with $\lambda=1$ is called a Steiner system, and the most celebrated ones 
have parameters $4$-$(11,5,1)$, $5$-$(12,6,1)$, $3$-$(22,6,1)$, $4$-$(23,7,1)$, and $5$-$(24,8,1)$, 
which are often referred to as the Mathieu designs or the Witt systems. 
In particular, the $4$-$(11,5,1)$ design $W_{11}$ arises from a 4-transitive action of the Mathieu group $M_{11}$ on 11 points, and its existence and uniqueness was first shown by Witt \cite{BJ}. The design $W_{11}$ has 66 blocks, and every two distinct blocks $B$ and $B'$ have $1$, $2$, or $3$ points in common. Let us define a binary symmetric relation $R_i$ on $W_{11}$ by 
\begin{equation*}
(B,B')\in R_i \Leftrightarrow \vert B\cap B'\vert=4-i,
\end{equation*}
for $i=1,2,3$ and let $R_0=\{(B,B)\mid B\in W_{11}\}$. 
Then the pair $\mathcal{W}=(W_{11},\{R_0,R_1,R_2,R_3\})$ is an association scheme (of 3 classes) by \cite[Theorem~5.25]{D}. Moreover, this scheme is $Q$-polynomial (see \cite{BI} for the definitions 
and more results about $P,Q$-polynomial association schemes). Williford \cite{GWW,W} compiled 
the tables of feasible parameters of primitive 3-class $Q$-polynomial association schemes 
(on up to 2800 vertices), where the uniqueness of the scheme $\mathcal{W}$ was left blank
(see also \cite[Appendix~B]{vD99}). 
In this note, we show that $\mathcal{W}$ is unique, i.e., it is determined up to isomorphism 
by its parameters. The proof is computer-assisted by Mathematica and relies on a spherical 
representation of the scheme \cite{BBB}.

\section{The parameters of $\mathcal{W}$}
Let us recall some standard facts from the theory of association schemes (see \cite{BI}).   
Let $A_i$ denote the logical matrix of the relation $R_i$, for $i=0,1,2,3$. Then $A_0$ is the identity matrix of size $66$ and:
\begin{enumerate}
\item $\sum_{i=0}^3 A_i = J$, the square all-one matrix of size $66$,
\item $A_i^\top=A_i$ ($0 \le i \le 3$),
\item $A_iA_j=\sum_{k=0}^3 p_{i,j}^kA_k$,
where $p_{i,j}^k$ are nonnegative integers  ($0 \le i,j,k \le 3$), 
called the {\bf intersection numbers} of the scheme, which we refer to as the parameters of the scheme.
\end{enumerate}

The intersection numbers of $\mathcal{W}$, written in the form of matrices $\left(L_i\right)_{kj}=(p_{i,j}^k)$,
are found in the tables by Williford \cite{W}:
\begin{align*}
L_1&=\left[
\begin{array}{cccc}
 0 & 30 & 0 & 0 \\
 1 & 15 & 10 & 4 \\
 0 & 15 & 6 & 9 \\
 0 & 8 & 12 & 10 \\
\end{array}
\right],\quad 
L_2=\left[
\begin{array}{cccc}
 0 & 0 & 20 & 0 \\
 0 & 10 & 4 & 6 \\
 1 & 6 & 10 & 3 \\
 0 & 12 & 4 & 4 \\
\end{array}
\right],\quad 
L_3=\left[
\begin{array}{cccc}
 0 & 0 & 0 & 15 \\
 0 & 4 & 6 & 5 \\
 0 & 9 & 3 & 3 \\
 1 & 10 & 4 & 0 \\
\end{array}
\right],
\end{align*}
and they determine (\cite[Theorem~4.1]{BI}) the first and second eigenmatrices $P$ and $Q=\frac{1}{66}P^{-1}$ of $\mathcal{W}$:
\begin{align*}
P&=\left[
\begin{array}{cccc}
 1 & 30 & 20 & 15 \\
 1 & 8 & -2 & -7 \\
 1 & -1 & -2 & 2 \\
 1 & -6 & 8 & -3 \\
\end{array}
\right],\quad
Q=\left[
\begin{array}{cccc}
 1 & 10 & 44 & 11 \\
 1 & \frac{8}{3} & -\frac{22}{15} & -\frac{11}{5} \\
 1 & -1 & -\frac{22}{5} & \frac{22}{5} \\
 1 & -\frac{14}{3} & \frac{88}{15} & -\frac{11}{5} \\
\end{array}
\right], 
\end{align*}
where the $(P)_{ij}$-entry ($0\le i,j\le 3$) is the eigenvalue of $A_j$ on the $i$-th maximal common eigenspace of the matrices $A_0,\ldots,A_3$ (which commute and hence can be simultaneously diagonalized). 

\begin{theorem}\label{thm:s1}
An association scheme with the above parameters is isomorphic to $\mathcal{W}$. 
\end{theorem}

In what follows, we assume that $\mathcal{X}=(X,\{R_0,R_1,R_2,R_3\})$ is an association scheme 
with the same parameters as $\mathcal{W}$. 
Let $E_1$ denote the orthogonal projection matrix onto the 1st maximal common eigenspace 
(of multiplicity $m_1=Q_{01}=10$) of the matrices $A_i$'s of $\mathcal{X}$.  
Since $E_1$ is positive semidefinite, 
we may regard $\frac{\vert X\vert}{m_1}E_1$ as the Gram matrix of vectors in the unit sphere $S^{m_1-1}$ in $\mathbb{R}^{m_1}$, and write $\frac{\vert X\vert }{m_1}E_1=F_1 F_1^\top$ where $F_1$ is a $\vert X\vert \times m_1$ matrix.  
We now identify a vertex $x\in X$ with the $x$-th row $\mathbf{x}\in \mathbb{R}^{m_1}$ of $F_1$, and such a map is said to be a {\bf spherical representation} 
of the scheme $\mathcal{X}$ into its first eigenspace.
Define the \textbf{angle set} 
$A(X):=\left\{\langle \mathbf{x},\mathbf{y}\rangle \mid {x},{y}\in X, {x}\neq {y}\right\}$ 
and observe that $A(X)=\left\{\frac{Q_{j1}}{Q_{01}} \mid 1\leq j\leq 3 \right\}$. 
Note that the map is injective, since $Q_{j1}\neq Q_{01}$ provided that $j\neq 0$; 
thus, it defines all relations of the scheme. 
Therefore, to prove Theorem \ref{thm:s1} it suffices to show that the Gram matrix 
of the vertex set $X$ embedded into the unit sphere $S^9$ in $\mathbb{R}^{10}$ 
is unique, up to orthogonal 
transformation.

\section{
Proof of Theorem \ref{thm:s1}
}

Observe that the graph $(X,R_2)$ is strongly regular with parameters $(66,20,10,4)$. 
Such a strongly regular graph is isomorphic to the triangular graph $T(12)$ by \cite{Chang}; 
thus, we may identify the point set $X$ with the vertex set of $T(12)$, 
which is ${[12]\choose 2}$, where $[12]:=\{1,\ldots, 12\}$, 
and $R_2$ with the edge set of $T(12)$, which is 
$\left\{\{v,v'\} \mid \vert v\cap v'\vert=1, v,v'\in {[12]\choose 2}\right\}$.
The next lemma can easily be seen from this description of $T(12)$.  

\begin{lemma}\label{lm:clique}
Each vertex of the triangular graph $T(12)$ is contained in two maximum cliques of order $11$. Each vertex outside such a clique has exactly two neighbors in it.
\end{lemma}

Note that a clique of order $11$ in $T(12)$ is a {\bf Delsarte clique}, 
as it attains the Delsarte bound \cite[Proposition 1.3.2]{BCN}. 
Fix a vertex $x=\{1,2\}$ of $T(12)$. The two Delsarte cliques containing $x$ are 
\begin{align*}
C_1=\{\{1,j\} \mid 2\leq j\leq 12\} \quad \text{ and } \quad C_2=\{\{2,j\} \mid j=1,3\leq j\leq 12\}. 
\end{align*}

Consider the image of $C_1$ in the spherical representation in $\mathbb{R}^{10}$: 
since $C_1$ is a clique, the angle between any two vectors is the same; 
thus, these 11 vectors form a regular simplex. 
Let $V_1$ be an $11\times 10$ matrix whose row vectors correspond 
to the vertices of $C_1$. Up to orthogonal transformation, we may assume that 
the matrix $V_1$ has the following form: 
\begin{align*}
V_1=\left[
\begin{smallmatrix}
 1 & 0 & 0 & 0 & 0 & 0 & 0 & 0 & 0 & 0 \\
 -\frac{1}{10} & \frac{3 \sqrt{11}}{10} & 0 & 0 & 0 & 0 & 0 & 0 & 0 & 0 \\
 -\frac{1}{10} & -\frac{\sqrt{11}}{30} & \frac{2 \sqrt{55}}{15} & 0 & 0 & 0 & 0 & 0 & 0 & 0 \\
 -\frac{1}{10} & -\frac{\sqrt{11}}{30} & -\frac{\sqrt{55}}{60} & \frac{\sqrt{385}}{20} & 0 & 0 & 0 & 0 & 0 & 0 \\
 -\frac{1}{10} & -\frac{\sqrt{11}}{30} & -\frac{\sqrt{55}}{60} & -\frac{\sqrt{385}}{140} & \frac{\sqrt{1155}}{35} & 0 & 0 & 0 & 0 & 0 \\
 -\frac{1}{10} & -\frac{\sqrt{11}}{30} & -\frac{\sqrt{55}}{60} & -\frac{\sqrt{385}}{140} & -\frac{\sqrt{1155}}{210} & \frac{\sqrt{33}}{6} & 0 & 0 & 0 & 0 \\
 -\frac{1}{10} & -\frac{\sqrt{11}}{30} & -\frac{\sqrt{55}}{60} & -\frac{\sqrt{385}}{140} & -\frac{\sqrt{1155}}{210} & -\frac{\sqrt{33}}{30} & \frac{\sqrt{22}}{5} & 0 & 0 & 0 \\
 -\frac{1}{10} & -\frac{\sqrt{11}}{30} & -\frac{\sqrt{55}}{60} & -\frac{\sqrt{385}}{140} & -\frac{\sqrt{1155}}{210} & -\frac{\sqrt{33}}{30} & -\frac{\sqrt{22}}{20} & \frac{\sqrt{330}}{20} & 0 & 0 \\
 -\frac{1}{10} & -\frac{\sqrt{11}}{30} & -\frac{\sqrt{55}}{60} & -\frac{\sqrt{385}}{140} & -\frac{\sqrt{1155}}{210} & -\frac{\sqrt{33}}{30} & -\frac{\sqrt{22}}{20} & -\frac{\sqrt{330}}{60} & \frac{\sqrt{165}}{15} & 0 \\
 -\frac{1}{10} & -\frac{\sqrt{11}}{30} & -\frac{\sqrt{55}}{60} & -\frac{\sqrt{385}}{140} & -\frac{\sqrt{1155}}{210} & -\frac{\sqrt{33}}{30} & -\frac{\sqrt{22}}{20} & -\frac{\sqrt{330}}{60} & -\frac{\sqrt{165}}{30} & \frac{\sqrt{55}}{10} \\
 -\frac{1}{10} & -\frac{\sqrt{11}}{30} & -\frac{\sqrt{55}}{60} & -\frac{\sqrt{385}}{140} & -\frac{\sqrt{1155}}{210} & -\frac{\sqrt{33}}{30} & -\frac{\sqrt{22}}{20} & -\frac{\sqrt{330}}{60} & -\frac{\sqrt{165}}{30} & -\frac{\sqrt{55}}{10} \\
\end{smallmatrix}
\right],
\end{align*} 
where the first row corresponds to $x$.

Now the problem is to add the remaining $55$ vectors such that together with the vectors from $C_1$ 
their Gram matrix is permutation equivalent to $E_1$. Since $A(X)=\{\alpha_1,\alpha_2,\alpha_3\}$, 
where $\alpha_1=\frac{4}{15},\alpha_2=-\frac{1}{10},\alpha_3=-\frac{7}{15}$, 
every such a vector $\mathbf{u}\in S^9$ in question satisfies 
$\displaystyle{V_1 \mathbf{u}^\top \in \{\alpha_1,\alpha_2,\alpha_3\}^{10}}$, 
and there are at most $3^{10}$ candidates for $\mathbf{u}$.

We first determine the coordinates of the vectors in the image of $C_2$. 
Let $z\in C_1 \setminus \{x\}$ and $y\in C_2$. 
Then $(x,z) \in R_2$, $(x,y)\in R_2$ and $(y,z)\in R_1\cup R_2 \cup R_3$. 
Since $p_{2,1}^2=6,p_{2,2}^2=10,p_{2,3}^2=3$ and, by Lemma \ref{lm:clique}, 
for each vertex $z\in C_1\setminus\{x\}$, there is exactly one vertex $y\in C_2\setminus\{x\}$ 
such that $\langle \mathbf{y},\mathbf{z}\rangle=\alpha_2$, it follows that
the vertices in $C_2$ are taken from the following set: 
\begin{align*}
Y_1&=\left\{\mathbf{y}\in S^9 \mid V_1 \mathbf{y}^\top =(\alpha_2,\mathbf{v})^\top,\mathbf{v}=\left(\{[\alpha_1]^6,[\alpha_2]^1,[\alpha_3]^3\}\right) \right\},   
\end{align*}
where $\left(\{[\alpha_1]^6,[\alpha_2]^1,[\alpha_3]^3\}\right)$ denotes a
vector of length 10 having $6$ entries equal to $\alpha_1$, $1$ entry $\alpha_2$, and $3$ entries $\alpha_3$. 
Note that $\vert Y_1\vert =840$. 
Consider a graph with vertex set $Y_1$ and edge set $E_1$ defined by $\{\mathbf{y},\mathbf{y}'\}\in E_1$ if and only if 
$\langle \mathbf{y},\mathbf{y}'\rangle=\alpha_2$. 
Then every clique of order $10$ in the graph $(Y_1,E_1)$ is a candidate for $C_2$.
Let us define a $10\times 10$ matrix $V_2=C\cdot V_1$, where  
\[
C=\frac{1}{3}
\begin{bmatrix*}[r]
0 & 1  & 1  & 1  & 1  & 1  & 1  & 0  & -1 & -1 & -1 \\
0 & 1  & 1  & 1  & 1  & -1 & -1 & -1 & 1  & 1  & 0 \\
0 & 1  & 1  & 0  & -1 & 1  & -1 & 1  & 1  & -1 & 1 \\
0 & 1  & 1  & -1 & 0  & -1 & 1  & 1  & -1 & 1  & 1 \\
0 & 1  & -1 & 1  & -1 & 1  & 1  & -1 & 0  & 1  & 1 \\
0 & 1  & -1 & -1 & 1  & 1  & 0  & 1  & 1  & 1  & -1 \\
0 & 0  & -1 & 1  & 1  & -1 & 1  & 1  & 1  & -1 & 1 \\
0 & -1 & 1  & 1  & -1 & 0  & 1  & 1  & 1  & 1  & -1 \\
0 & -1 & 1  & -1 & 1  & 1  & 1  & -1 & 1  & 0  & 1 \\
0 & -1 & 0  & 1  & 1  & 1  & -1 & 1  & -1 & 1  & 1 \\
\end{bmatrix*},
\]
and $V=\left[\begin{matrix}
 V_1\\
 V_2
\end{matrix}\right]$. We can prove the following lemma with the aid of a computer.
\begin{lemma}\label{lem:q1}
There exist exactly $30240$ cliques of order $10$ in the graph $(Y_1,E_1)$. 
For every such clique $C_2$, the Gram matrix of $C_1\cup C_2$ is permutationally 
equivalent to that of the rows of $V$.
\end{lemma}

By Lemma \ref{lem:q1}, we may fix the set $Y_2$ of the row vectors of $V_2$ representing $C_2$ and extend the above argument to other Delsarte cliques, which yields that $X\setminus (C_1 \cup C_2)$ can be represented by a subset of the following set: 
\begin{align*}
Y=\left\{\mathbf{u}\in S^9 \mid V_1 \mathbf{u}^\top =\mathbf{v}^\top,\mathbf{v}=\left(\{[\alpha_1]^6,[\alpha_2]^2,[\alpha_3]^3\}\right)\right\}\setminus Y_2.
\end{align*} 

Note that $\vert Y\vert =4610$, which is somewhat smaller than $3^{10}$. We determine the set $Y$ and  select only those of its vectors whose inner products with the 
vectors from $Y_2$ belong to $A(X)$. 
Let $Z$ denote the set of such vectors. 
With the aid of a computer, we obtain $\vert Z\vert =90$. 
We proceed by finding a maximal subset of $Z$ such that the inner product of every two of its distinct vectors is in $A(X)$. 
\begin{lemma}
The set $Z$ splits into $Z_1\cup Z_2$ in such a way that $\vert Z_1\vert =\vert Z_2\vert =45$, $Z_1\cap Z_2=\emptyset$, $A(Z_1)=A(Z_2)=\{\alpha_1,\alpha_2,\alpha_3\}$, and $\langle \mathbf{z}_1,\mathbf{z}_2\rangle \not\in\{\alpha_1,\alpha_2,\alpha_3\}$ for all $\mathbf{z}_1\in Z_1$ and all $\mathbf{z}_2\in Z_2$. Moreover, the Gram matrices $C_1\cup C_2 \cup Z_1$ and $C_1\cup C_2 \cup Z_2$ are permutationally equivalent.  
\end{lemma}

Thus, we may assume that the $66$ vertices of $X$ are represented by the row vectors from $V_1$, $V_2$, and $Z_1$, and we can directly verify that, together with the binary relations defined by their inner products, they form an association scheme of $3$ classes with the same parameters as $\mathcal{W}$. This completes the proof of Theorem~\ref{thm:s1}.

\subsubsection*{Acknowledgments}
The authors would like to thank the reviewers for valuable comments. 
The research of Alexander Gavrilyuk is supported by JSPS KAKENHI Grant Number 22K03403.
The research of Sho Suda is supported by JSPS KAKENHI Grant Number 22K03410.



\end{document}